\documentclass[letterpaper, 10 pt, conference]{ieeeconf}

\IEEEoverridecommandlockouts
\usepackage{amsmath,amssymb,bm}

\usepackage{amsthm}
\usepackage{graphicx,setspace}
\usepackage[colorlinks=true,allcolors=blue]{hyperref}
\usepackage{tikz}
\usetikzlibrary{arrows.meta,positioning,calc,fit,backgrounds}

\newtheoremstyle{ieeerunin}{3pt}{3pt}{\itshape}{}{\itshape}{:}{.5em}{}
\theoremstyle{ieeerunin}

\newtheorem{definition}{Definition}
\newtheorem{lemma}{Lemma}
\newtheorem{theorem}{Theorem}
\newtheorem{proposition}{Proposition}
\newtheorem{corollary}{Corollary}

\newenvironment{pf}%
  {\par\addvspace{3pt}\noindent\hspace*{\parindent}\textit{Proof:}\ \ignorespaces}%
  {\hfill$\blacksquare$\par\addvspace{3pt}}

\newcommand{\Real}{\mathbb{R}}
\newcommand{\Zset}{\mathcal{Z}}
\newcommand{\dd}{\mathrm{d}}

\newcommand{\relu}[1]{\left[#1\right]_{+}}
\newcommand{\xu}{x_{u}}
\newcommand{\xsl}{x_{1}}
\newcommand{\xsr}{x_{2}}

\title{\LARGE Translate, then decide: asymptotic optimality of \\ latent-space control in bistable systems}

\author{Arthur Castello B. de Oliveira$^{1}$
\thanks{$^{1}$Center for Computational Neuroscience at the Flatiron Institute.\ \tt\small acastello@flatironinstitute.org}%
}

\begin{document}
\maketitle
\thispagestyle{empty}
\pagestyle{empty}

\begin{abstract}
This paper asks whether an optimal decision can be separated from the physical dynamics on which it acts. It studies this question through minimum-energy barrier crossing in scalar bistable systems. It characterizes the finite-horizon optimal feedback and shows that, as the horizon grows, it converges to a law that doubles and reverses the drift while it opposes crossing, and switches off once the barrier is passed. In physical coordinates, this law inherits the plant's nonlinearity, entangling decision and dynamics. We show that a plant-specific translation of state, input, and time maps every double well in the considered class to the same unit-rate linear world model, allowing a downstream decision stage to use a single plant-independent control primitive. This separation generally loses optimality at finite horizons because the translation does not preserve the physical control energy. Nevertheless, we prove that the excess physical energy of the layered controller vanishes as the horizon grows and that its limiting feedback reproduces the asymptotic optimum exactly. Thus, translation and optimization generally do not commute at a finite horizon, but commute in the asymptotic limit.
\end{abstract}

\section{Introduction}\label{sec:intro}

A growing view in neural computation treats neural activity as the state of a population evolving through a high-dimensional space. In the prefrontal cortex (PFC), nonlinear mixed selectivity yields high-dimensional representations accessible to simple readouts~\cite{rigotti2013importance}, while population activity during context-dependent decisions follows structured dynamical trajectories~\cite{mante2013context}. More broadly, high-dimensional neural activity can encode low-dimensional task structure~\cite{vyas2020computation,chung2021neural}, and learning can reshape representational geometry so that common structure is exposed across related tasks~\cite{bernardi2020geometry,wakhloo2026neural}.

Evidence from learning and generalization suggests that such a structure can be reused across different physical realizations of a task. Prior experience facilitates learning of structurally related problems~\cite{harlow1949formation,braun2009motor}, while neural representations can become increasingly abstract with learning~\cite{tse2007schemas,zhou2021evolving}. Particularly relevant here, mice trained on reversal-learning problems sharing an abstract structure but differing in their physical implementation exhibit representations in medial PFC that generalize across problems, while hippocampal representations retain more problem-specific information~\cite{samborska2022complementary}. Computational accounts similarly separate slowly acquired reusable structure from fast inference within a particular task~\cite{wang2018prefrontal,whittington2020tolman}. Together, these observations motivate a layered abstraction: a task-general component reused across problems, and a situation-specific representation that relates it to the physical problem at hand.

If neural computation is viewed as controlled population dynamics, control theory provides a natural language for asking what such a separation costs. Optimal feedback control has long provided models of sensorimotor and neural computation~\cite{todorov2002optimal,madhav2020synergy}, including recently at the level of individual neurons~\cite{moore2024neuron}. From the control perspective, abstraction is closely related to the use of coordinate and feedback transformations to simplify nonlinear dynamics~\cite{isidori1995nonlinear}. Such transformations, however, also modify the input channel and therefore the control metric. As such, simplifying the dynamics does not, in general, preserve the associated optimal-control problem. The relationship between feedback linearization and quadratic optimal control has been studied explicitly, including asymptotic regimes in which a feedback-linearizing controller approaches an optimal controller~\cite{braslavsky1998near}. More generally, exact preservation of optimality under composition has been established for particular classes of problems~\cite{todorov2009compositionality}. These observations motivate the question of when a reusable transformed representation can support, rather than merely approximate, optimal decision-making.

We study this question in the simplest setting where representation, decision, and optimality can all be treated analytically: minimum-energy barrier crossing in a scalar bistable system $\dot x = f(x)+u$
Two stable equilibria are separated by an unstable equilibrium $x_u$, and the task is to cross this barrier using minimum control energy. Minimum-energy steering between coexisting attractors and across basin boundaries has been studied in nonlinear dynamical systems~\cite{silchenko2004energy}, while bistable dynamics have long served as models of commitment between competing alternatives~\cite{wang2002probabilistic}. The scalar setting considered here isolates the transverse crossing problem from the additional geometric choices available along a higher-dimensional separatrix.

Minimum-energy transfer across unstable equilibria was recently studied in~\cite{behera2026nonlinear}. There, optimal state-feedback laws were derived near one- and two-dimensional linearized saddles and shown to reduce, in the free-horizon limit, to one-sided linear primitives. An appendix further derives the corresponding scalar nonlinear result, establishing the form of optimal nonlinear scalar control. The resulting primitive, however, remains tied to the particular dynamics and transfer specification, and reducing the degree of task conditioning is identified as an open direction. Here we take that nonlinear scalar result as a starting point and ask a different question: \emph{can plant-specific information be confined to a representation layer, leaving the decision itself unchanged across different nonlinear plants?}

We first formulate the nonlinear problem directly in terms of the barrier, and characterize the finite-horizon minimum-energy crossing law and its asymptotic limit, recovering the ``double-or-nothing’’ structure of~\cite{behera2026nonlinear} -- reverse the drift with twice its strength while it opposes crossing, and apply no control once it assists it. 

Our main result is that this strategy can be broken into two stages (translation and decision), but that doing so has a nontrivial optimality cost. We show that the resulting layered controller is generally suboptimal over finite horizons and characterize exactly when equality is possible. Nevertheless, as the crossing horizon tends to infinity, we observe that the discrepancy vanishes and the universal decision law recovers the asymptotic optimum exactly. Thus, for this class of problems, translation and optimization generally do not commute at a finite horizon, but commute exactly in the asymptotic limit.

\section{Theoretical framework}\label{sec:framework}

\subsection{Bistable scalar systems}\label{sec:dw-def}

We study the scalar control-affine plant
\begin{equation}\label{eq:plant}
  \dot x=f(x)+u,\qquad x\in\Real,\quad u\in\Real ,
\end{equation}
with $f\in C^{1}(\Real;\Real)$, and write $x(\cdot\,;x_{0},u)$ for the solution from $x_{0}$ under the input $u$. We specifically focus on scalar ``double wells'' defined as:

\begin{definition}[Double well]\label{def:dw}
  The plant \eqref{eq:plant} is a \emph{double well} on an open interval $I\subseteq\Real$ if $f$ has exactly three zeros in $I$, two of which are asymptotically stable and one unstable for the uncontrolled flow $\dot x=f(x)$, with the unstable being hyperbolic. We denote the stable zeros by $\xsl,\xsr$ and the unstable one by $\xu$.
\end{definition}

The canonical instance is the cubic $f(x)=x-x^{3}$, with $\xsl=-1$, $\xsr=1$ and $\xu=0$, but every statement below is proved for Definition~\ref{def:dw} without further assumptions required. First, from Definition~\ref{def:dw}, the following result can easily be proven by continuity and the intermediate value theorem.

\begin{proposition}[Structure of a double well]\label{prop:structure}
  Let \eqref{eq:plant} be a double well on $I\subset\mathbb{R}$. Then, $\xsl<\xu<\xsr$, and $f$ is sign definite on each of the intervals $(\xsl,\xu)$ and $(\xu,\xsr)$, with
  \begin{equation}\label{eq:signstruct}
    f(x)<0 \ \text{ on } (\xsl,\xu),
    \qquad
    f(x)>0 \ \text{ on } (\xu,\xsr) .
  \end{equation}
  Consequently the basins of attraction of $\xsl$ and $\xsr$ are $\{x\in I~|~x<\xu\}$ and $\{x\in I~|~x>\xu\}$, separated by the single point $\xu$.
\end{proposition}


Uncontrolled, the state is therefore carried to the bottom of whichever well it starts in, and the problem of interest is to move it to the other well, using as little control energy as possible, that is, to \emph{cross the barrier}. Before formally writing the problem, we next consider the barrier closely and rewrite the system in a form convenient for expressing the barrier crossing problem.

\subsection{The barrier coordinate and the transport rate}\label{sec:dw-barrier}

The unstable equilibrium separates the two basins, and it is convenient to measure the state by its signed distance from it. Define the \emph{barrier coordinate} $\varphi(x):=x-\xu$, whose zero set $\Zset:=\{\xu\}$ we call the \emph{barrier}. Since $f(\xu)=0$, the barrier is invariant under the drift, thus by continuity no uncontrolled trajectory crosses it, and the sign of $\varphi$ records which basin the state occupies. In higher dimensions, the barrier will be a hypersurface -- the separatrix between the two basins -- and $\varphi$ a function vanishing on it, but in one dimension, both objects are as simple as they can be.

Because $f$ vanishes at $\xu$, the quotient
\begin{equation}\label{eq:a}
  a(x):=\frac{f(x)}{x-\xu}\quad(x\neq\xu),
  \qquad a(\xu):=f'(\xu),
\end{equation}
is well defined and continuous, since the value at $\xu$ being the limit $f(x)/(x-\xu)\to f'(\xu)$ is guaranteed by differentiability. We call $a$ the \emph{transport rate} and might omit both its and $\varphi$'s dependency on $x$ if clear by context. Furthermore, we might write $a(t)$ and $\varphi(t)$ if  we want to emphasize their time evolution, but they should be understood as $a(x(t))$ and $\varphi(x(t))$. 

By construction, the transport rate satisfies the \emph{transport identity}
\begin{equation}\label{eq:transport}
  f(x)=a(x)\,\varphi(x)\qquad\text{for all }x ,
\end{equation}
so that along the uncontrolled flow
\begin{equation}\label{eq:free}
  \dot\varphi=a(x)\,\varphi ,
\end{equation}
and at the barrier
\begin{equation}\label{eq:lambda}
  a(\xu)=f'(\xu)=:\lambda>0 .
\end{equation}
For the cubic, $a(x)=1-x^{2}$ and $\lambda=1$.

From \eqref{eq:free} we can formalize the notion that the barrier is \emph{repelling} as follows:
\begin{corollary}\label{cor:repel}
  On $(\xsl,\xsr)$ one has $a(x)>0$, and consequently
  \begin{equation}\label{eq:repel}
    \frac{\dd}{\dd t}\varphi^{2}=2a\varphi^{2}>0
    \qquad\text{for } x\in(\xsl,\xsr)\setminus\{\xu\}
  \end{equation}
  along the free drift. Moreover $a(\xsl)=a(\xsr)=0$.
\end{corollary}


In other words, $|\varphi|$ grows under the free drift on \emph{both} sides of the barrier, and control action is required to change the sign of $\varphi$.  Also note that the repulsion necessarily weakens near the bottom of the wells, where $a\to0$. As such, we fix an open interval of interest for the crossing problem.
\begin{equation}\label{eq:U}
  \xu\in U\subset(\xsl,\xsr),\qquad \overline{U}\subset(\xsl,\xsr),
\end{equation}
on which continuity of $a$ and Corollary~\ref{cor:repel} give
\begin{equation}\label{eq:abounds}
  0<a_{-}:=\min_{x\in\overline{U}}a(x) \le a(x) \le a_{+}:=\max_{x\in\overline{U}}a(x)<\infty .
\end{equation}

Finally, also note that in the controlled plant the input enters the barrier coordinate directly, that is, differentiating $\varphi(x)$ along \eqref{eq:plant} and using
\eqref{eq:transport} gives
\begin{equation}\label{eq:channel}
  \dot\varphi=a(x)\,\varphi+u,
\end{equation}
and each value of $\dot\varphi$ therefore determines $u$ point-wise uniquely, and conversely.

\section{The barrier-crossing problem}\label{sec:ocp}

In this section, we formally introduce the barrier-crossing problem for a double well, discuss well-posedness issues, and provide analytic solutions for different versions of the problem.

\subsection{Problem statement}\label{sec:ocp-def}

Fix $x_{0}\in U$ with $\varphi_{0}:=\varphi(x_{0})>0$ (so the state starts in the positive basin) and a horizon $\tau>0$. Let $\mathcal{U}_{\tau}(x_{0})$ be the set of inputs $u\in L^{2}([0,\tau];\Real)$ whose resulting trajectory remains in $\overline{U}$ throughout $[0,\tau]$ -- i.e. $x(t,x_0,u)\in \overline U,~\forall t\in[0,\tau]$. The \emph{$\tau$-crossing problem} is defined as
\begin{equation}\label{eq:Vtau}
  \begin{aligned}
    V_{\tau}(x_{0}):=\inf_{u\in\mathcal{U}_{\tau}(x_{0})}
      &\quad J_{\tau}(u):=\tfrac12\int_{0}^{\tau}u(t)^{2}\,\dd t ,\\
    \text{subject to}&\quad \varphi\bigl(x(\tau;x_{0},u)\bigr)<0 ,
  \end{aligned}
\end{equation}
with the infimum over an empty feasible set equal to $+\infty$.

Notice that, for the problem definition in \eqref{eq:Vtau}, enlarging the horizon can only help -- given $\tau_1<\tau_2$ all feasible $\tau_1$-crossing solutions are also feasible for the $\tau_2$-crossing -- so $\tau\mapsto V_{\tau}$ is non-increasing and nonnegative, hence convergent. This allows the definition of the \emph{asymptotic crossing value}
\begin{equation}\label{eq:V}
  V(x_{0}):=\lim_{\tau\to\infty}V_{\tau}(x_{0})=\inf_{\tau>0}V_{\tau}(x_{0}) .
\end{equation}

Also note that \eqref{eq:Vtau} is ill-posed in the sense that the infimum might exist but will never be attained, since the constraint set is open. We will consider this issue next.

\subsection{Attainment}\label{sec:ocp-attain}

Notice in problem \eqref{eq:Vtau} that crossing further past the barrier costs strictly more, so the cheapest admissible terminal states are those just barely across, and a minimizing sequence of feasible inputs drives $\varphi(x(\tau))$ to zero from below without any member achieving the infimum. This observation is made explicit through a well-posed surrogate, obtained by requiring a margin $\varepsilon>0$, with $\xu-\varepsilon\in U$, replacing the endpoint conditions by
\begin{equation}\label{eq:eps}
  \varphi(x_{0})>0,
  \qquad \varphi\bigl(x(\tau;x_{0},u)\bigr)\le-\varepsilon ,
\end{equation}
and writing $V^{\varepsilon}_{\tau}(x_{0})$ for the resulting value. For this related $\varepsilon,\tau$-crossing problem, attainment can be proven through the following lemma.

\begin{lemma}[Attainment with a margin]\label{lem:attain}
  Fix $\tau>0$ and $\varepsilon>0$. Then, the infimum defining $V^{\varepsilon}_{\tau}(x_{0})$ is attained. Moreover $\varepsilon\mapsto V^{\varepsilon}_{\tau}(x_{0})$ is non-decreasing and $V^{\varepsilon}_{\tau}(x_{0})\to V_{\tau}(x_{0})$ as $\varepsilon\to0^{+}$.
\end{lemma}

\begin{pf}
\emph{(sketch.)}
Let $(u_n)$ be a minimizing sequence. Since $J_\tau(u_n)\to V^\varepsilon_\tau(x_0)<\infty$, the sequence is bounded in $L^2([0,\tau])$. By reflexivity, a subsequence converges weakly to some $u^\star\in L^2([0,\tau])$. Let $x_n:=x(\cdot\,;x_0,u_n)$. Admissibility gives $x_n(t)\in\overline U$, while boundedness of $f$ on $\overline U$ and the $L^2$ bound on $u_n$ imply equicontinuity by Cauchy--Schwarz. Hence, by the Arzel\`a--Ascoli theorem, a further subsequence converges uniformly to some $x^\star:[0,\tau]\to\overline U$. Passing to the limit in the integral dynamics, using uniform convergence of $f(x_n)$ and weak convergence of $u_n$, shows that $x^\star=x(\cdot\,;x_0,u^\star)$. The terminal constraint is preserved by continuity of $\varphi$, so $u^\star$ is admissible. Finally, weak lower semicontinuity of the $L^2$ norm gives $J_\tau(u^\star)\le\liminf_{n\to\infty}J_\tau(u_n)=V^\varepsilon_\tau(x_0)$, while admissibility gives the reverse inequality. Thus the infimum is attained.

For the limit as $\varepsilon\to0^+$, the admissible sets increase as $\varepsilon$ decreases and remain contained in the admissible set of \eqref{eq:Vtau}, so $V_\tau(x_0)\le V^\varepsilon_\tau(x_0)$. Conversely, for any $\delta>0$, choose $u$ feasible for \eqref{eq:Vtau} with $J_\tau(u)<V_\tau(x_0)+\delta$. Since its terminal state satisfies $\varphi(x(\tau;x_0,u))<0$, the same $u$ is feasible for the $\varepsilon,\tau$-problem for all sufficiently small $\varepsilon>0$. Hence $V_\tau(x_0)\le V^\varepsilon_\tau(x_0)\le J_\tau(u)<V_\tau(x_0)+\delta$. Since $\delta>0$ is arbitrary, $V^\varepsilon_\tau(x_0)\to V_\tau(x_0)$ as $\varepsilon\to0^+$.
\end{pf}

Thus, for each $\varepsilon>0$ there is a genuine optimal control that crosses strictly, with value converging to $V_{\tau}$ from above. This observation enables us to look at the explicit solutions of the crossing problem through the Pontryagin maximum principle (PMP).

\subsection{Closed-form solution}\label{sec:ocp-solve}

We will next look at the exact solution of the $\varepsilon,\tau$-crossing problem and analyze how this solution behaves as we approach the $\tau$-crossing problem and the asymptotic crossing problem (both of which are unattainable). Throughout, $x_{0}\in U$ with $\varphi_{0}:=\varphi(x_{0})>0$, and we regard $a$ as a function of the barrier coordinate, writing $a(\varphi)$ for $a(\xu+\varphi)$.

\begin{theorem}[The $\varepsilon,\tau$-optimal feedback]\label{thm:mu}
  For every $\tau>0$ there exists a $\mu>0$ for which the stationary feedback
  \begin{equation}\label{eq:umu}
    u_{\mu}(x)\;=\;-a(x)\varphi(x)-\sqrt{a(x)^{2}\varphi(x)^{2}+2\mu},
  \end{equation}
  is the unique solution of the $\varepsilon,\tau$-crossing problem.
  Furthermore, the closed loop is given by
  \begin{equation}\label{eq:wmu}
    \dot\varphi=-\sqrt{a(\varphi)^{2}\varphi^{2}+2\mu}=:-w_{\mu}(\varphi)\;<\;0 ,
  \end{equation}
  and attains the following optimal value:
  \begin{equation}\label{eq:Jmu}
    V^{\varepsilon}_{\tau}(x_{0})
    =\int_{-\varepsilon}^{\varphi_{0}}
       \Bigl[w_{\mu}(\varphi)+a(\varphi)\varphi\Bigr]\dd\varphi\;-\;\mu\,\tau .
  \end{equation}
\end{theorem}

\begin{pf}
First notice that any minimizing trajectory must be monotone. Indeed, a reversal in one dimension produces a repeated state, and deleting the intervening loop strictly decreases the energy -- the lost time can be restored at no cost by holding the state at \(x_u\), where \(f(x_u)=0\), before resuming the trajectory. Hence a minimizer remains in \([x_u-\varepsilon,x_0]\subset U\), so the state constraint is inactive. By Lemma~\ref{lem:attain} a minimizer exists, to which the Pontryagin maximum principle applies. Let \(\nu\ge0\) denote the multiplier of the terminal constraint \(\varphi(\tau)+\varepsilon\le0\). If the constraint were inactive, complementary slackness would give \(\nu=0\), hence \(p(\tau)=0\); since \(\dot p=-pf’(x)\), this would imply \(p\equiv0\), and stationarity would give \(u\equiv0\), which cannot cross the invariant barrier. Thus the terminal constraint binds, \(\varphi(\tau)=-\varepsilon\).

With Hamiltonian $H=p\,(f+u)+\tfrac12u^{2}$, first order stationarity in $u$ gives $u=-p$ (a minimum, as $\partial_{u}^{2}H=1>0$), and the costate obeys $\dot p=-pf'(x)$. Substituting $u=-p$ gives $H=pf-\tfrac12p^{2}=:-\mu$, which is constant along the extremal because the problem is autonomous and the horizon is fixed. Solving $p^{2}-2pf-2\mu=0$ gives $p=f\pm\sqrt{f^{2}+2\mu}$, hence $u=-f\mp\sqrt{f^{2}+2\mu}$ and $\dot\varphi=f+u=\mp\sqrt{f^{2}+2\mu}$. Since this is the only feasible solution to also satisfy the PMP, it is the unique minimizer.

Next we show that $\mu>0$ through the Hamiltonian. Along the extremal, $-\mu=H=pf-\tfrac12p^{2}$ with $u=-p$, so $\mu=uf+\tfrac12u^{2}=u\bigl(f+\tfrac{u}{2}\bigr)$; substituting $u=\dot\varphi-f$ and $f = a\varphi$ gives $\mu=\tfrac12\bigl(\dot\varphi^{2}-f^{2}\bigr)=\tfrac12\bigl(\dot\varphi^{2}-a^{2}\varphi^{2}\bigr)$,
{constant}. Since $\varphi$ is continuous with $\varphi(0)=\varphi_{0}>0$ and $\varphi(\tau)=-\varepsilon<0$, there is $t^{*}\in(0,\tau)$ with $\varphi(t^{*})=0$, where $f=a\varphi$ vanishes. At that time we have $\mu=\tfrac12\dot\varphi(t^{*})^{2}>0$. Therefore $f^{2}+2\mu>0$ everywhere, so the two branches $p=f\pm\sqrt{f^{2}+2\mu}$ never meet and no switching between them is possible. In other words, $\dot\varphi=\mp\sqrt{f^{2}+2\mu}$ has constant sign, negative by the boundary conditions, which selects the lower branch. Substituting $f=a\varphi$ recovers \eqref{eq:umu} and \eqref{eq:wmu}.

To compute the optimal value, reparametrize time by $\varphi$ through $\dd t=-\dd\varphi/w_{\mu}$, so that with $u=\dot\varphi-a\varphi$ the cost is
\[
  J=\int_{-\varepsilon}^{\varphi_{0}}
    \frac{\bigl(w_{\mu}+a\varphi\bigr)^{2}}{2\,w_{\mu}}\,\dd\varphi .
\]
Expanding the integrand and using $a^{2}\varphi^{2}=w_{\mu}^{2}-2\mu$ gives $\frac{(w_{\mu}+a\varphi)^{2}}{2w_{\mu}}=w_{\mu}+a\varphi-\frac{\mu}{w_{\mu}}$,
and $\int_{-\varepsilon}^{\varphi_{0}}\dd\varphi/w_{\mu}=\tau$ by \eqref{eq:wmu}, giving \eqref{eq:Jmu}.
\end{pf}

Although Theorem \ref{thm:mu} gives a closed form solution for the $\varepsilon,\tau$-crossing problem, the exact meaning behind the undefined constant $\mu>0$ is unclear. We clarify it in the next result.

\begin{theorem}[The horizon--multiplier relation]\label{thm:mutau}
  Fix $\varepsilon>0$ and $\xu-\varepsilon\in U$ with $\varphi_0>0$. The multiplier $\mu$ of Theorem~\ref{thm:mu} and the horizon $\tau$ are related by
  \begin{equation}\label{eq:Tmu}
    \tau=T(\mu):=\int_{-\varepsilon}^{\varphi_{0}}
      \frac{\dd\varphi}{\sqrt{a(\varphi)^{2}\varphi^{2}+2\mu}} ,
  \end{equation}
  where $T:(0,\infty)\to(0,\infty)$ is continuous, strictly decreasing and bijective. In particular
  \begin{equation}\label{eq:Tlimits}
    \lim_{\mu\to0^{+}}T(\mu)=\infty ,
    \quad
    \lim_{\mu\to\infty}T(\mu)=0 ,
    \quad
    \lim_{\mu\to0^{+}}\mu\,T(\mu)=0 ,
  \end{equation}
  so every horizon $\tau>0$ is realized by exactly one $\mu=T^{-1}(\tau)>0$, and $T^{-1}$ is continuous and strictly decreasing.
\end{theorem}

\begin{pf}
Equation \eqref{eq:Tmu} is the statement that the crossing takes time $\tau$: by \eqref{eq:wmu} the closed loop obeys $\dot\varphi=-w_{\mu}(\varphi)$ with $w_{\mu}>0$, so $\varphi$ is a strictly decreasing, hence invertible, function of time, and $\dd t=-\dd\varphi/w_{\mu}(\varphi)$ integrates over $\varphi\in[-\varepsilon,\varphi_{0}]$ to the elapsed time, recovering \eqref{eq:Tmu}.

Next, notice that for $\mu>0$ the integrand is continuous and bounded above by $(2\mu)^{-1/2}$ on the compact interval $[-\varepsilon,\varphi_{0}]$, so $T(\mu)$ is finite; and it is bounded below by $(a_{+}^{2}\Lambda^{2}+2\mu)^{-1/2}(\varphi_{0}+\varepsilon)>0$ with $\Lambda:=\max\{\varphi_{0},\varepsilon\}$, using \eqref{eq:abounds}. Hence $T$ maps $(0,\infty)$ into $(0,\infty)$.

Next notice that for each fixed $\varphi$ the map $\mu\mapsto(a^{2}\varphi^{2}+2\mu)^{-1/2}$ is continuous and strictly decreasing in $\mu$ pointwise. Hence, $T$ is continuous and strictly decreasing by dominated convergence on compact subsets of $(0,\infty)$.

\emph{The limit as $\mu\to\infty$.} Bounding the integrand above by $(2\mu)^{-1/2}$,
\[
  0<T(\mu)\le\frac{\varphi_{0}+\varepsilon}{\sqrt{2\mu}}\ \longrightarrow\ 0 .
\]

{Next we look at $\mu\to0^{+}$.} Let $(\mu_{n})$ be any strictly decreasing sequence with $\mu_{n}\to0^{+}$, and write $h_{n}(\varphi):=\bigl(a(\varphi)^{2}\varphi^{2}+2\mu_{n}\bigr)^{-1/2}$. Since $\mu_{n}$ decreases, the sequence $h_{n}$ is nonnegative and nondecreasing in $n$, and converges pointwise on $[-\varepsilon,\varphi_{0}]\setminus\{0\}$ to $h(\varphi):=1/\bigl(a(\varphi)|\varphi|\bigr)$. By the monotone convergence theorem,
\[
  T(\mu_{n})=\int_{-\varepsilon}^{\varphi_{0}}h_{n}
  \ \longrightarrow\ \int_{-\varepsilon}^{\varphi_{0}}h
  =\int_{-\varepsilon}^{\varphi_{0}}\frac{\dd\varphi}{a(\varphi)|\varphi|}=\infty.
\]
Hence $T(\mu_{n})\to\infty$, and since $T$ is monotone the full limit $\lim_{\mu\to0^{+}}T(\mu)=\infty$ follows.

These estimates for the limits also control the product of the multiplier with the horizon. Bounding the integrand of \eqref{eq:Tmu} above by $(2\mu)^{-1/2}$ gives $T(\mu)\le(\varphi_{0}+\varepsilon)/\sqrt{2\mu}$ for every $\mu>0$, whence
\begin{equation}\label{eq:muT}
  0\ \le\ \mu\,T(\mu)\ \le\ \frac{(\varphi_{0}+\varepsilon)\sqrt{\mu}}{\sqrt2}
  \ \xrightarrow[\ \mu\to0^{+}\ ]{}\ 0 .
\end{equation}
Thus, although $T(\mu)$ diverges as $\mu\to0^{+}$, it does so slowly enough that $\mu\,T(\mu)$ vanishes.

Finally, since $T$ is continuous and strictly decreasing on the interval $(0,\infty)$, it is injective, and by \eqref{eq:Tlimits} and the intermediate value theorem its range is all of $(0,\infty)$. It is therefore a bijection onto $(0,\infty)$, and its inverse is continuous and strictly decreasing.
\end{pf}


\begin{theorem}[The asymptotic limit]\label{thm:limit}
  Fix $\varepsilon>0$ with $\xu-\varepsilon\in U$, and let $x_{0}\in U$ with
  $\varphi_{0}>0$. As $\tau\to\infty$ the multiplier satisfies
  $\mu=T^{-1}(\tau)\to0^{+}$, the feedback \eqref{eq:umu} converges uniformly
  on $\overline{U}$ to
  \begin{equation}\label{eq:u0}
    \begin{aligned}
        u_{0}(x)&=-a(x)\varphi(x)-a(x)\bigl|\varphi(x)\bigr|
        \\&=\begin{cases}
           -2a(x)\varphi(x)=-2f(x), & \varphi(x)>0,\\[2pt]
           \ \ 0, & \varphi(x)<0,
         \end{cases}
    \end{aligned}
  \end{equation}
  and the value converges to
  \begin{equation}\label{eq:V0}
  \footnotesize
    \begin{aligned}
        V^{\varepsilon}(x_{0})&=\lim_{\tau\to\infty}V^{\varepsilon}_{\tau}(x_{0})
        =\int_{0}^{\varphi_{0}}2\,a(\varphi)\,\varphi\,\dd\varphi
        =2\int_{\xu}^{x_{0}}f(\sigma)\,\dd\sigma.
    \end{aligned}
  \end{equation}
  In particular, both the limiting law and the limiting value are independent
  of $\varepsilon$.
\end{theorem}

\begin{pf}
By Theorem~\ref{thm:mutau} the map $T$ is a decreasing bijection of $(0,\infty)$ onto itself, so its inverse is decreasing with $T^{-1}(\tau)\to0^{+}$ as $\tau\to\infty$. This proves the first claim, and it allows us to replace the limit $\tau\to\infty$ by the limit $\mu\to0^{+}$ throughout.

Next, for $\alpha,\beta\ge0$ one has $\sqrt{\alpha+\beta}\le\sqrt{\alpha}+\sqrt{\beta}$ by subadditivity of the square root, so with $\alpha=a^{2}\varphi^{2}$ and $\beta=2\mu$,
\begin{equation}\label{eq:wunif}
  0\ \le\ w_{\mu}(\varphi)-a(\varphi)|\varphi|
    \ =\ \sqrt{a^{2}\varphi^{2}+2\mu}-\sqrt{a^{2}\varphi^{2}}
    \ \le\ \sqrt{2\mu},
\end{equation}
the left inequality holding because $\mu>0$. The bound is uniform in $\varphi$, so $w_{\mu}\to a|\varphi|$ uniformly on $\overline{U}$ as $\mu\to0^{+}$, and therefore $u_{\mu}=-a\varphi-w_{\mu}$ converges uniformly to $-a\varphi-a|\varphi|$, which is the first line of \eqref{eq:u0}. The two cases follow from the two signs of $\varphi$, using $a>0$ on $U$ by Corollary~\ref{cor:repel}: for $\varphi>0$ the terms add to $-2a\varphi=-2f$, and for $\varphi<0$ they cancel exactly. 

{Next, for evaluating the convergence of the value function,} consider each term in \eqref{eq:Jmu} separately. For the integral, \eqref{eq:wunif} gives
\begin{equation*}
  \left|\int_{-\varepsilon}^{\varphi_{0}}\bigl[w_{\mu}+a\varphi\bigr]\dd\varphi
   -\int_{-\varepsilon}^{\varphi_{0}}\bigl[a|\varphi|+a\varphi\bigr]\dd\varphi
  \right|
  \le (\varphi_{0}+\varepsilon)\sqrt{2\mu},
\end{equation*}
which goes to zero with $\mu$. The limiting integrand $a|\varphi|+a\varphi$ vanishes identically for $\varphi<0$ and equals $2a\varphi$ for $\varphi>0$, so the integral converges to $\int_{0}^{\varphi_{0}}2a\varphi\,\dd\varphi$. For the remaining term, $\mu\tau=\mu\,T(\mu)\to0$ by \eqref{eq:Tlimits}. Combining the two limits gives the second line of \eqref{eq:V0}, and the third follows from $a\varphi=f$ and the change of variable $\varphi=\sigma-\xu$.

Finally, neither limit retains $\varepsilon$, since the law \eqref{eq:u0} never involved it, and in the value the segment $[-\varepsilon,0]$ contributes nothing because the limiting integrand vanishes there.
\end{pf}

We emphasize that both limits in Theorem~\ref{thm:limit} are independent of $\varepsilon$. The reason is visible in the integrand of \eqref{eq:Jmu}: at $\mu=0$ the contribution of the segment $[-\varepsilon,0]$ is $a|\varphi|+a\varphi=0$, so the far side of the barrier costs nothing, the repelling drift carrying the state onward by itself. Accordingly $V^{\varepsilon}(x_{0})=V(x_{0})$ for every $\varepsilon>0$, and the asymptotic $\varepsilon$-crossing problem, the asymptotic crossing problem and the barrier-\emph{reaching} problem with target $\Zset$ all share the value \eqref{eq:V0} and the limiting law \eqref{eq:u0}, with the caveat that only the latter is actually attainable.

In other words, they do not share a solution. For every $\mu>0$ the feedback \eqref{eq:umu} crosses in finite time and is optimal for its horizon, but the limit \eqref{eq:u0} does not cross at all, since on the near side it renders $\dot\varphi=-a\varphi$, which reaches the barrier only as $t\to\infty$. The asymptotic value is thus approached but never attained, and $u_{0}$ is best read as the optimal \emph{reaching} law, of which the crossing controls \eqref{eq:umu} are finite-horizon perturbations of size $\sqrt{2\mu}$ in the crossing rate at the barrier.


\subsection{The optimal law entangles plant and decision}
\label{sec:entanglement}

For a \emph{linear} unstable plant $\dot z=\lambda z+u$ with $\lambda>0$, the same problem has optimal solution $u=-2\lambda\relu{z}$~\cite[Sec.~II]{behera2026nonlinear}. This is the special case of \eqref{eq:u0} with constant transport rate $a\equiv\lambda$, and the plant enters the law only through the scalar gain $\lambda$.

For a nonlinear plant the dependence is far richer. Where $\varphi>0$ the law is $u_{0}=-2f(x)$ (for the cubic, $-2(x-x^{3})$), so it contains the entire drift and not a single parameter. The crossing task is the same for every double well, yet the law computed for one double well does not carry over to another: changing the plant changes the whole function to be implemented. In physical coordinates, then, identifying the plant and optimizing the decision are entangled.

The next section shows that this conclusion is an artifact of the coordinates in which \eqref{eq:u0} is written, and can be circumvented through a two-layer interpretation.

\section{Translation and decision}\label{sec:twolayer}

The optimal law \eqref{eq:u0} entangles knowledge of the plant with the decision. This section asks whether the two can be separated, as in the layered organization motivated in the introduction, and, if so, what the separation costs in terms of optimality are.

\subsection{An alternative: translate first, then decide}\label{sec:translate}

The difficulty with \eqref{eq:u0} is that a single map must know both \emph{where} the barrier lies and \emph{how strongly} the drift opposes the crossing. Suppose instead that these are supplied by an upstream stage, and that the controller downstream of it solves a problem posed entirely in the barrier coordinate.

Concretely, let the command be issued as a multiple of the local transport rate,
\begin{equation}\label{eq:decomp}
  u=a(x)\,d ,
\end{equation}
with $d$ a abstract \emph{decision} variable. Substituting into \eqref{eq:channel} we get $\dot\varphi=a(x)\bigl(\varphi+d\bigr)$, so that $a>0$ multiplies both the state and the command and can be removed by a change of clock. Define the \emph{transport time}
\begin{equation}\label{eq:clock}
  s(t):=\int_{0}^{t}a\bigl(x(\sigma)\bigr)\,\dd\sigma ,
\end{equation}
strictly increasing since $a\ge a_{-}>0$ on $\overline{U}$ by \eqref{eq:abounds}. In this clock $\varphi$ evolves as
\begin{equation}\label{eq:lti}
  \frac{\dd\varphi}{\dd s}=\varphi+d ,
\end{equation}
a linear time-invariant unstable plant with unit rate and unit input gain -- and, crucially, \emph{the same plant for every double well}. All information about the particular drift has been absorbed into the translation \eqref{eq:decomp} and the clock \eqref{eq:clock}, and the linear system \eqref{eq:lti} acts as a ``world model'' for this type of task.

\begin{figure}[t]
\centering
\begin{tikzpicture}[
  >={Stealth[length=2mm]},
  font=\small,
  box/.style  = {draw, rounded corners=1.5pt, minimum height=9mm,
                 minimum width=16mm, inner xsep=4pt, align=center},
  gain/.style = {draw, circle, inner sep=1pt, minimum size=6mm},
  lbl/.style  = {font=\scriptsize, inner sep=1.5pt},
  scale=0.8,
  transform shape
]
\node[box, fill=black!5] (plant) {plant\\$\dot x=f(x)+u$};
\node[box, right=16mm of plant]      (phi)  {$\varphi(x)$};
\node[box, below=5mm of phi]         (rate) {$a(x)$};
\node[box, right=16mm of phi, fill=black!12] (dec) {$-2\relu{\,\cdot\,}$};
\node[gain, below=6.5mm of dec]        (mul)  {$\times$};

\begin{scope}[on background layer]
  \node[draw, dashed, rounded corners=3pt, fill=black!2,
        fit=(phi)(rate), inner sep=6pt] (trans) {};
\end{scope}

\draw[->] (plant.east) -- node[lbl,above]{$x$} (phi.west);
\draw[->] (phi.east)   -- node[lbl,above]{$\varphi$} (dec.west);
\coordinate (tap) at ($(plant.east)!0.45!(phi.west)$);
\draw[->] (tap) |- (rate.west);
\fill (tap) circle (1pt);
\draw[->] (dec.south) -- node[lbl,right]{$d$} (mul.north);
\draw[->] (rate.east) -- node[lbl,above]{$a$}  (mul.west);
\coordinate (uin) at ($(plant.south)+(0,-18mm)$);
\draw[->] (mul.south) |- (uin) -- node[lbl,right,pos=0.6]{$u=a\,d$}
          (plant.south);

\node[lbl, anchor=south] at (trans.north)
      {translation: $\varphi$ and $a$};
\node[lbl, anchor=south] at (dec.north)   {decision};
\end{tikzpicture}
\caption{Translation and decision. The translation stage supplies the barrier coordinate $\varphi$ and the transport rate $a$, both learned from the free dynamics; the decision stage is the context-agnostic solution of the barrier crossing for system \eqref{eq:lti}, and identical for every double well. Scaling the abstract decision by $a$ produces the command $u=a\,d$, which is equivalent to running the decision in the transport clock \eqref{eq:clock}, in which the loop reduces to the unit-rate plant \eqref{eq:lti} regardless of $f$.\\ \vspace{-20pt}}
\label{fig:twolayer}
\end{figure}
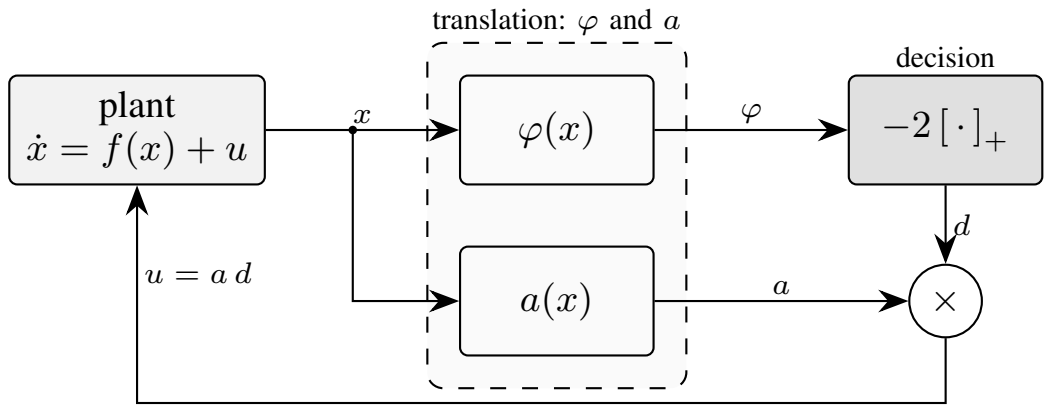

This suggests the two-stage architecture of Fig.~\ref{fig:twolayer}, which we call \emph{translation and decision}. The translation stage represents $\varphi$ and $a$ and performs \eqref{eq:decomp}. The decision stage sees only the channel \eqref{eq:lti} and its own clock, and solves the crossing problem there, in the variable available to it and with the cost natural to it:
\begin{equation}\label{eq:virtual}
  \min_{d}\ \tfrac12\int_{0}^{S}d(s)^{2}\,\dd s
  \quad\text{s.t.}\quad \eqref{eq:lti},\ \
  \varphi(S)=-\varepsilon .
\end{equation}
We emphasize that \eqref{eq:virtual} is what a downstream stage can actually pose. It has no access to $a$, that being the property of the plant the translation stage was introduced to absorb, and it measures duration in the clock it runs on. Whether this is the right objective is precisely the question of this section.

The convenience of a linear representation, however, lies exactly on the simplicity of its solution. Problem \eqref{eq:virtual} is the scalar unstable minimum-energy problem of \cite[Sec.~II]{behera2026nonlinear} with unit rate and unit input gain, whose solution is available in closed form. Specifically, for horizon $S$ and margin $\varepsilon$, the optimal decision for $s\in[0,S]$ is
\begin{equation}\label{eq:gramian}
  d^{\star}_S(s)=-B\,e^{-s},\quad B:=\frac{\varepsilon+e^{S}\varphi_{0}}{\sinh S}>0,
\end{equation}
with virtual cost $\tfrac12\int_0^S d^{\star2}_S\,\dd s=(\varepsilon+e^{S}\varphi_{0})^{2}/(e^{2S}-1)$. This follows from a least-norm problem. 
As $S\to\infty$ and $\varepsilon\to0^{+}$, $B\to2\varphi_0$, so $d^{\star}_S(s)\to-2\varphi_{0}e^{-s}$. Along the limiting closed loop $\dd\varphi/\dd s=\varphi+d=-\varphi$ one has $\varphi(s)=\varphi_{0}e^{-s}$, so the limit is the stationary feedback
\begin{equation}\label{eq:dstar}
  d^{\star}=-2\varphi ,
\end{equation}
which reflects the open-loop pole $+1$ of \eqref{eq:lti} to $-1$. Since the state is already across the barrier where $\varphi<0$, the feedback acts only on the positive side, so \eqref{eq:dstar} reads $d^{\star}=-2\relu{\varphi}$: a rectified linear law whose gain is exactly $2$. Because \eqref{eq:lti} is the same world model for every double well, this decision is identical for all plants, and all plant-specific information remains in the translation.

\subsection{On the equivalence between approaches}\label{sec:neq}

One can easily see that the decision stage optimizes the wrong quantity. By \eqref{eq:decomp} and $\dd s=a\,\dd t$, the true control energy is not \eqref{eq:virtual} but a \emph{weighted} energy in the same variable and clock,
\begin{equation}\label{eq:weights}
    \small
  \tfrac12\int u^{2}\,\dd t
  \;=\;\tfrac12\int d^{2}a^{2}\,\dd t
  \;=\;\tfrac12\int d^{2}\,a\,\dd s \;\neq\;\tfrac12\int d^2\dd s,
\end{equation}
the weight being exactly the transport rate that the change of clock was designed to remove. This is the energy premium usually charged against feedback linearization: a virtual constraint is imposed on a chosen output, the drift's contribution is cancelled in its entirety, and the desired motion re-supplied from scratch; unless the re-supply is priced with the true weights, unnecessary energy is injected into the system, resulting in suboptimality.


\begin{theorem}[Translation and optimization do not commute]\label{thm:neq}
  Fix $\varepsilon>0$ with $\xu-\varepsilon\in U$ and $x_0\in U$ with $\varphi_0:=\varphi(x_0)>0$. For $S>0$, let $d^\star_S$ be the solution \eqref{eq:gramian} of the decision problem \eqref{eq:virtual} over the transport horizon $S$, and let $x^{\mathrm{ts}}_S$ be the trajectory of the layered closed loop
  \[
    \dot x=f(x)+a(x)\,d^\star_S(s),\quad \dot s=a(x),\quad x(0)=x_0,\ \ s(0)=0,
  \]
  with layered law $u^{\mathrm{ts}}_S(t):=a(x^{\mathrm{ts}}_S(t))\,d^\star_S(s(t))$. Let $\tau_S$ be the first time at which $\varphi(x^{\mathrm{ts}}_S(t))=-\varepsilon$, and let $u_\mu$ be the optimal $\varepsilon,\tau_S$-crossing feedback \eqref{eq:umu} of Theorem~\ref{thm:mu}, with $\mu=T^{-1}(\tau_S)>0$ by Theorem~\ref{thm:mutau}. Then
  \begin{enumerate}
    \item $\tau_S$ is finite, and $u^{\mathrm{ts}}_S$ restricted to $[0,\tau_S]$ is a feasible solution of the $\varepsilon,\tau_S$-crossing problem;
    \item $u^{\mathrm{ts}}_S$ is optimal for the $\varepsilon,\tau_S$-crossing problem if and only if $a(x^{\mathrm{ts}}_S(t))$ is constant on $[0,\tau_S]$;
    \item For every horizon $\tau>0$ there exists $S>0$ with $\tau=\tau_S$, and any such $S$ satisfies $S\ge a_-\tau$; in particular, $S_\tau\to\infty$ as $\tau\to\infty$ for any choice of $S_\tau$ with $\tau_{S_\tau}=\tau$.
  \end{enumerate}
\end{theorem}

\begin{pf}
Let $\varphi_S(s):=e^{s}\varphi_0-B\sinh s$, the unique solution of $\dd\varphi/\dd s=\varphi+d^{\star}_S(s)$, $\varphi(0)=\varphi_0$
with $\varphi_S(S)=-\varepsilon$. Computing $\mu$ for the LTI model \eqref{eq:lti} gives
\begin{equation}\label{eq:conserved}
  \mu_{\textrm{LTI}}=\tfrac12 d^{\star}_S(s)^{2}+\varphi_S(s)\,d^{\star}_S(s)=\tfrac12B(B-2\varphi_0)\quad,
\end{equation}
for all $s\in[0,S]$. Simplifying the expression gives,
\[
  \mu_{\textrm{LTI}}=\frac12\frac{(\varepsilon+e^{S}\varphi_0)(\varepsilon+e^{-S}\varphi_0)}{\sinh^{2}S}>0 .
\]
By \eqref{eq:conserved}, $(\varphi_S+d^{\star}_S)^{2}=\varphi_S^{2}+2\mu_{\textrm{LTI}}>0$, so $\varphi_S+d^{\star}_S$ has no zero on $[0,S]$. Since $\varphi_S(S)= -\varepsilon<0<\varphi_0=\varphi_S(0)$, by the mean value theorem we conclude
\begin{equation}\label{eq:negative}
  \tfrac{\dd}{\dd s}\varphi_S=\varphi_S+d^{\star}_S<0\ \text{ on }[0,S],
\end{equation}
so $\varphi_S$ is strictly decreasing from $\varphi_0$ to $-\varepsilon$, and $\varphi_S([0,S])=[-\varepsilon,\varphi_0]$.

Next, since $\xu-\varepsilon, x_0\in U$, then $K:=[\xu-\varepsilon,x_0]\subset U$, and $\xu+\varphi_S(s)\in K$ for all $s\in[0,S]$. On $K$, $a\ge a_->0$ by \eqref{eq:abounds}, so
$t(s):=\int_0^{s}({a\bigl(\xu+\varphi_S(\sigma)\bigr)})^{-1}{\dd\sigma},~ s\in[0,S]$,
is $C^1$ and strictly increasing. Let $s(t)$ be its inverse on $[0,t(S)]$ and set $x^{\mathrm{ts}}_S(t):=\xu+\varphi_S(s(t))$. Then $\dot s=1/t'(s)=a(x^{\mathrm{ts}}_S)$ and, using $f=a\varphi$ and $\varphi_S'=\varphi_S+d^\star_S$,
\[
  \dot x^{\mathrm{ts}}_S=\varphi_S'(s)\,\dot s=a\,(\varphi_S+d^\star_S)=f(x^{\mathrm{ts}}_S)+a(x^{\mathrm{ts}}_S)\,d^\star_S(s),
\]
so $(x^{\mathrm{ts}}_S,s)$ solves the layered loop on $[0,t(S)]$ and stays in $K\subset U$. Likewise, uniqueness follows from standard arguments and is omitted for space reasons.

Next we look at the crossing time. The loop is run until $s$ reaches $S$, at real time $t(S)$. Since $\varphi_S$ is strictly decreasing and $\varphi_S(s)=-\varepsilon$ only at $s=S$, the first time at which $\varphi=-\varepsilon$ is
\begin{equation}\label{eq:tauS}
  \tau_S=t(S)=\int_0^{S}\frac{\dd\sigma}{a\bigl(\xu+\varphi_S(\sigma)\bigr)}\in\Bigl[\tfrac{S}{a_+},\tfrac{S}{a_-}\Bigr],
\end{equation}
which is finite. We now proceed to proving each claim (not in order).

\emph{Proof of 1).} We established that $x^{\mathrm{ts}}_S(t)\in[\xu-\varepsilon,x_0]\subset\overline U$ for $t\in[0,\tau_S]$, $\varphi_0>0$, and $\varphi\bigl(x^{\mathrm{ts}}_S(\tau_S)\bigr)=-\varepsilon$. The law $u^{\mathrm{ts}}_S$ is continuous and bounded on $[0,\tau_S]$, hence in $L^2$. Since $x^{\mathrm{ts}}_S$ solves $\dot x=f(x)+u^{\mathrm{ts}}_S(t)$ (uniquely), $x^{\mathrm{ts}}_S=x(\cdot\,;x_0,u^{\mathrm{ts}}_S)$. Hence $u^{\mathrm{ts}}_S$ is feasible for the $\varepsilon,\tau_S$-crossing problem.

\emph{Proof of 3).} Substituting $\sigma=S\theta$ in \eqref{eq:tauS},
\[
  \tau_S=S\int_0^1\frac{\dd\theta}{a\bigl(\xu+\varphi_S(S\theta)\bigr)} .
\]
The map $(S,\theta)\mapsto\varphi_S(S\theta)\in[-\varepsilon,\varphi_0]$, is continuous, and the integrand is continuous and bounded between $1/a_+$ and $1/a_-$, therefore $S\mapsto\tau_S$ is continuous on $(0,\infty)$. By \eqref{eq:tauS}, $\tau_S\to0$ as $S\to0^+$ and $\tau_S\to\infty$ as $S\to\infty$, so by the intermediate value theorem every $\tau>0$ equals $\tau_S$ for some $S$. If $\tau_S=\tau$, then \eqref{eq:tauS} gives $\tau\le S/a_-$, i.e., $S\ge a_-\tau$, so $S_\tau\to\infty$ as $\tau\to\infty$ for any choice of $S_\tau$.

\emph{Proof of 2). ($\Rightarrow$).} By Lemma~\ref{lem:attain} a minimizer of the $\varepsilon,\tau_S$-crossing problem exists, and by the proof of Theorem~\ref{thm:mu} every minimizer is the signal generated by the closed loop of $u_\mu$ for some $\mu>0$, with $\mu=T^{-1}(\tau_S)$ fixed by Theorem~\ref{thm:mutau}; so this signal $u^\star$ is the unique minimizer. Suppose now that $u^{\mathrm{ts}}_S$ is optimal. Then $u^{\mathrm{ts}}_S=u^\star$ almost everywhere, hence everywhere by continuity. By uniqueness of the open-loop trajectory, $x^{\mathrm{ts}}_S=x_\mu$, so $u^{\mathrm{ts}}_S(t)=u_\mu(x^{\mathrm{ts}}_S(t))$. Dividing by $a>0$ and using \eqref{eq:umu},
\[
  d^{\star}_S(s)=-\varphi-\sqrt{\varphi^{2}+2\mu/a^{2}}\;\Rightarrow\; d^{\star2}_S+2\varphi\,d^{\star}_S=\frac{2\mu}{a^{2}} .
\]
By \eqref{eq:conserved} the left-hand side equals $2\mu_{\textrm{LTI}}$ for all $t$, hence $a\bigl(x^{\mathrm{ts}}_S(t)\bigr)^{2}=\mu/\mu_{\textrm{LTI}}$ is constant on $[0,\tau_S]$.

\emph{Proof of 2). ($\Leftarrow$).} Suppose $a\bigl(x^{\mathrm{ts}}_S(t)\bigr)\equiv a_c$ and set $\tilde\mu:=a_c^{2}\mu_{\textrm{LTI}}>0$. By \eqref{eq:conserved}, $(d^\star_S+\varphi)^2=\varphi^2+2\mu_{\textrm{LTI}}$, and by \eqref{eq:negative} $d^\star_S+\varphi<0$, so
\[
  d^{\star}_S=-\varphi-\sqrt{\varphi^{2}+2\mu_{\textrm{LTI}}}=-\varphi-\sqrt{\varphi^{2}+2\tilde\mu/a_c^{2}} .
\]
Multiplying by $a_c=a(x^{\mathrm{ts}}_S(t))$,
\[
  u^{\mathrm{ts}}_S(t)=-a\varphi-\sqrt{a^{2}\varphi^{2}+2\tilde\mu}=u_{\tilde\mu}\bigl(x^{\mathrm{ts}}_S(t)\bigr).
\]
Hence $x^{\mathrm{ts}}_S$ solves $\dot x=f(x)+u_{\tilde\mu}(x)$, $x(0)=x_0$, and by uniqueness it is the closed-loop trajectory of $u_{\tilde\mu}$. Along it $\dot\varphi=-\sqrt{f^2+2\tilde\mu}<0$, so the first time at which $\varphi=-\varepsilon$ is $T(\tilde\mu)$ by \eqref{eq:Tmu}. This first-hitting time is $\tau_S$ by definition, and $\tau_S=T(\mu)$ by the choice of $\mu$ in the statement. Hence $T(\tilde\mu)=T(\mu)$, and since $T$ is injective (Theorem~\ref{thm:mutau}), $\tilde\mu=\mu$. Therefore $x^{\mathrm{ts}}_S=x_\mu$ and $u^{\mathrm{ts}}_S=u^\star$, so $u^{\mathrm{ts}}_S$ is optimal.

Finally, since $u^{\mathrm{ts}}_S$ is feasible (by 1)), $J_{\tau_S}(u^{\mathrm{ts}}_S)\ge V^{\varepsilon}_{\tau_S}(x_0)$, and since the optimal solution is unique, the inequality is strict if and only if $a\bigl(x^{\mathrm{ts}}_S(t)\bigr)$ is not constant.
\end{pf}

Notice that the true optimum in the virtual variable is $d_{\mu}=-\varphi-\sqrt{\varphi^{2}+2\mu/a^{2}}$, in which the transport rate enters through the term $2\mu/a^{2}$. However, as $\mu$ is the Hamiltonian constant of Theorem~\ref{thm:mu}, it goes to zero as the horizon grows -- $\mu=T^{-1}(\tau)\to0^{+}$ by Theorem~\ref{thm:mutau} -- and the term vanishes. This intuition is formalized in the following result.

\begin{theorem}[Translation and optimization commute asymptotically]
\label{thm:factor}
  Fix $\varepsilon>0$ with $\xu-\varepsilon\in U$ and $x_0\in U$ with $\varphi_0:=\varphi(x_0)>0$. For each $S>0$, let $u^{\mathrm{ts}}_S$ and $\tau_S$ be the layered control and its crossing time defined in Theorem~\ref{thm:neq}. Then, as $S\to\infty$, $\tau_S\to\infty$ and
  \begin{equation}\label{eq:costgap}
    J_{\tau_S}\!\left(u^{\mathrm{ts}}_S\right)
    -V^{\varepsilon}_{\tau_S}(x_0)
    \longrightarrow 0 .
  \end{equation}
  Moreover, the corresponding limiting decision is $d^\star=-2\relu{\varphi}$, and its composition with the translation \eqref{eq:decomp} reproduces the asymptotic reaching law of Theorem~\ref{thm:limit}:
  \begin{equation}\label{eq:compose}
    a(x)d^\star
    =-2a(x)\relu{\varphi(x)}
    =u_0(x),
    \qquad x\in\overline U .
  \end{equation}
\end{theorem}

\begin{pf}
By Theorem~\ref{thm:neq}, $\tau_S\to\infty$ as $S\to\infty$. Let $\nu:=\mu_{\mathrm{LTI}}(S)$ and $r_\nu(\varphi):=\sqrt{\varphi^2+2\nu}$. From \eqref{eq:gramian} and the relation $\dd\varphi/\dd s=\varphi+d_S^\star$, the layered trajectory in transport coordinates satisfies
\begin{equation*}
d_S^\star=-\varphi-r_\nu,\qquad
u_S^{\mathrm{ts}}=-a(\varphi)(\varphi+r_\nu),\qquad
\dot\varphi=-a(\varphi)r_\nu.
\end{equation*}
Changing variables from $t$ to $\varphi$ therefore gives
\begin{equation*}
J_{\tau_S}\!\left(u_S^{\mathrm{ts}}\right)
=
\int_{-\varepsilon}^{\varphi_0}
a(\varphi)
\frac{(\varphi+r_\nu(\varphi))^2}{2r_\nu(\varphi)}
\,\dd\varphi.
\end{equation*}
By Theorem~\ref{thm:limit},
\begin{equation*}
V(x_0)
=
\int_{-\varepsilon}^{\varphi_0}
a(\varphi)\bigl(|\varphi|+\varphi\bigr)\,\dd\varphi,
\end{equation*}
and hence
\begin{equation*}
J_{\tau_S}\!\left(u_S^{\mathrm{ts}}\right)-V(x_0)
=
\int_{-\varepsilon}^{\varphi_0}
a(\varphi)
\frac{\bigl(r_\nu(\varphi)-|\varphi|\bigr)^2}
     {2r_\nu(\varphi)}
\,\dd\varphi.
\end{equation*}
Since $a(\varphi)\le a_+$,
\begin{equation*}
\footnotesize
0\le
J_{\tau_S}\!\left(u_S^{\mathrm{ts}}\right)-V(x_0)
\le
a_+\int_{-\infty}^{\infty}
\frac{\bigl(\sqrt{\varphi^2+2\nu}-|\varphi|\bigr)^2}
     {2\sqrt{\varphi^2+2\nu}}
\,\dd\varphi
=a_+\nu,
\end{equation*}
where the last equality follows by symmetry and the substitution $\varphi=\sqrt{2\nu}\sinh y$ on $\varphi\ge0$. By \eqref{eq:conserved}, $\nu=\mu_{\mathrm{LTI}}(S)\to0$ as $S\to\infty$, so
\begin{equation*}
J_{\tau_S}\!\left(u_S^{\mathrm{ts}}\right)\to V(x_0).
\end{equation*}
Moreover, $u_S^{\mathrm{ts}}$ is feasible for the $\varepsilon,\tau_S$-crossing problem by Theorem~\ref{thm:neq}, while Theorem~\ref{thm:limit} gives
\begin{equation*}
V(x_0)\le V^\varepsilon_{\tau_S}(x_0)\le J_{\tau_S}\!\left(u_S^{\mathrm{ts}}\right).
\end{equation*}
Consequently,
\begin{equation*}
0\le
J_{\tau_S}\!\left(u_S^{\mathrm{ts}}\right)-V^\varepsilon_{\tau_S}(x_0)
\le
J_{\tau_S}\!\left(u_S^{\mathrm{ts}}\right)-V(x_0)
\to0,
\end{equation*}
which proves \eqref{eq:costgap}.

Finally, \eqref{eq:gramian} gives
\begin{equation*}
d_S^\star(\varphi)
=
-\varphi-\sqrt{\varphi^2+2\nu},
\end{equation*}
and therefore
\begin{equation*}
0\le
d_S^\star(\varphi)+2\relu{\varphi}
\le\sqrt{2\nu}
\end{equation*}
for all $\varphi$. Hence $d_S^\star\to-2\relu{\varphi}$ uniformly. Thus the limiting decision is $d^\star=-2\relu{\varphi}$, and $a(x)d^\star = -2a(x)\relu{\varphi(x)} = u_0(x)$ by Theorem~\ref{thm:limit}.
\end{pf}

Perhaps surprisingly, Theorem~\ref{thm:factor} shows that the layered strategy of translating the dynamics before optimizing still approaches an optimal strategy asymptotically. The decision stage optimizes an objective that is demonstrably generally not the true one (Theorem~\ref{thm:neq}) and yet in the limit it approaches the optimum exactly.

The practical content of this is visible in the relation between the two laws. By \eqref{eq:wunif} the finite-horizon optimum differs from the asymptotic one by at most a positive constant $\sqrt{2\mu}$. Furthermore, by Theorem~\ref{thm:mutau} the relationship $\tau=T(\mu)$ is a decreasing bijection of $(0,\infty)$ onto itself, and any perturbation of this kind, however small, converts the asymptotic reaching law into a finite-time crossing law, shortening the horizon monotonically as the perturbation grows.

Taken together, these two observations suggest a division of labor. The decision is made using a linear, context-independent world model, while everything plant-specific is confined to the translation layer, which supplies $\varphi$ and $a$. The learning problem is thereby displaced: rather than optimizing a control law in an infinite-dimensional function space as the plant is identified, one learns a representation and lets a fixed downstream primitive act on it. This analysis provides theoretical evidence for the value of such abstraction when studying interactions between different parts of the brain.

\section{Conclusions and future work}\label{sec:concl}

This paper provides theoretical evidence that a layered organization, in which a fixed decision stage acts on a world model that hides the specifics of the plant, can recover optimal control solutions for the asymptotic minimum-energy barrier crossing in generic scalar double wells.

We first solved the finite-time crossing problem with a margin. For every horizon, the minimum-energy solution is a stationary feedback with a single scalar parameter, in one-to-one correspondence with the horizon. As the horizon is extended, this feedback converges uniformly to a simple law: reverse the drift with twice its strength while it opposes the crossing, and switch off once the barrier is passed. Surprisingly, the limiting cost does not depend on the margin, and the limit is attained only as a reaching law, not as a crossing. In physical coordinates, this law inherits the plant's nonlinearity, entangling knowledge of the plant with the decision.

We then showed that the entanglement can be undone. Rewriting the input as a multiple of the local rate of instability and reparameterizing time in terms of that rate reduces every double well to the same linear time-invariant world model, whose minimum-energy solution is known in closed form~\cite{behera2026nonlinear}. Everything plant-specific is confined to the translation stage, while the decision stage is identical across all plants and is never relearned.

The decision stage, however, minimizes an energy that generally differs from the true one, and the resulting layered strategy is strictly suboptimal for the original problem except if the transport rate is constant. This loss is borne entirely by a horizon-dependent parameter, and vanishes as the horizon is relaxed to infinity. This means that the layered strategy coincides exactly with the global optimum in the asymptotic limit, in which case blindness to the specific plant costs nothing. This gives theoretical support to the intuition that what must be learned is the representation, not the optimal decision.

Finally, the construction is not specific to one dimension -- the barrier coordinate, the local rate of instability and the change of time are well defined for any plant whose uncontrolled dynamics have an invariant hypersurface separating two basins, and they reduce the crossing problem to the same world model. What does not carry over is global optimality: in higher dimensions, the cost of crossing can vary along the barrier, so the optimal controller may first move along it, and the residual suboptimality of the layered strategy depends on how the barrier is represented. Characterizing this dependence, and the case of barriers at hyperbolic saddles, is left to future work.

\bibliographystyle{plain}
\begin{spacing}{0.85}
\bibliography{bibliography.bib}
\end{spacing}

\end{document}